\documentclass{article}

\usepackage{a4,latexsym, amsfonts}
\usepackage{amssymb, amsmath}
\usepackage{amsfonts, latexsym}
\usepackage{booktabs}

\newcommand{\qed}{\hfill $\square$ \vspace{3mm}}
\newcommand{\proof}{\noindent \textbf{Proof \ }}
\newcommand{\im}{\mathrm{im\, }}
\newcommand{\rank}{\mathrm{rank\, }}

\newcommand{\fix}{\mathrm{fix\, }}
\newcommand{\dom}{\mathrm{dom\, }}
\newcommand{\opd}{\mathrm{opd\, }}
\newcommand{\ord}{\mathrm{ord\, }}

\newcommand{\Ker}{\mathrm{Ker\, }}

\newtheorem{theorem}{Theorem}
\newtheorem{lemma}[theorem]{Lemma}
\newtheorem{corollary}[theorem]{Corollary}

\newtheorem{proposition}[theorem]{Proposition}

\begin{document}
	
	\title{On certain semigroups of finite oriented and order-decreasing partial transformations}
	
	\author{Gonca Ay\i k$^{1}$,  Hayrullah Ay\i k$^{1}$, Ilinka Dimitrova$^{2,*}$, J\"{o}rg Koppitz$^{3}$ \vspace{3mm}\\
		$^{1}$Department of Mathematics, \c{C}ukurova University, Adana, Turkey \\
		$^{2}$Faculty of Mathematics and Natural Science, \\
		South-West University "Neofit Rilski", Blagoevgrad, Bulgaria \\
		$^{3}$Institute of Mathematics and Informatics,\\
		Bulgarian Academy of Sciences, Sofia, Bulgaria.}

	\date{}
	\maketitle

\let\thefootnote\relax\footnote{$^{*}$ Corresponding author.}

\let\thefootnote\relax\footnote{\, \ E-mail adresses:  agonca@cu.edu.tr (Gonca Ay\i k), hayik@cu.edu.tr (Hayrullah Ay\i k), ilinka\_dimitrova@swu.bg (Ilinka Dimitrova), koppitz@math.bas.bg (J\"{o}rg Koppitz)}	

\let\thefootnote\relax\footnote{\, \ The first and second authors were supported by Scientific and Technological Research Council of Turkey (TUBITAK) under the Grant Number 123F463. The first and second authors thank to TUBITAK for their supports.}

\begin{abstract}
	\noindent
	Let $\mathcal{PORD}_{n}$ be the semigroup consisting of all oriented and order-decreasing partial transformations
	on the finite chain $X_{n}=\{ 1<\cdots<n \}$. Let $\mathcal{IORD}_{n}$ be the subsemigroup of $\mathcal{PORD}_{n}$ consisting of all injective transformations of $\mathcal{PORD}_{n}$.
	For $2\leq r\leq n$, let $\mathcal{PORD}(n,r) =\{ \alpha\in \mathcal{PORD}_{n} :\lvert \im(\alpha) \rvert \leq r\}$ and $\mathcal{IORD}(n,r)=\{ \alpha \in \mathcal{IORD}_{n} :\lvert \im(\alpha )\rvert \leq r\}$. In this paper, we determine some minimal generating sets and ranks of $\mathcal{PORD}(n,r)$ and $\mathcal{IORD}(n,r)$, and moreover, we characterize the maximal subsemigroups of $\mathcal{PORD}(n,r)$ and $\mathcal{IORD}(n,r)$.
	
\end{abstract}

\vspace{3mm}

\vspace{3mm}


\section{Introduction}

For $n\in \mathbb{N}$, let $\mathcal{PT}_{n}$, $\mathcal{I}_{n}$, and $\mathcal{T}_{n}$ be the partial, injective, and full transformation semigroups on the finite chain $X_{n}=\{1,\ldots ,n\}$ under its natural order, respectively. An element $\alpha\in \mathcal{PT}_{n}$ with the domain $\dom (\alpha) =\{x_{1}< \cdots < x_{t}\}$ ($0\leq t\leq n$) (note that $\alpha$ is the empty transformation $0_n$ if $t=0$) is called \emph{order-preserving} (\emph{order-reversing}) if $x_{1}\alpha \leq\cdots\leq x_{t}\alpha$ ($x_{t}\alpha \leq\cdots\leq x_{1}\alpha$), and \emph{orientation-preserving} (\emph{orientation-reversing}) if the sequence $(x_{1}\alpha ,\ldots, x_{t}\alpha)$ is cyclic (anti-cyclic), i.e. there exists no more than one subscript $i$ such that $x_{i+1}\alpha <x_{i}\alpha$ ($x_{i}\alpha <x_{i+1}\alpha$), where $x_{t+1} = x_{1}$. Furthermore, $\alpha\in \mathcal{PT}_{n}$ is called \emph{order-decreasing} (\emph{extensive} or \emph{order-increasing}) if $x\alpha \leq x$ $(x\leq x\alpha)$ for all $x\in \dom(\alpha)$, \emph{monotone} if $\alpha$ is order-preserving or order-reversing, and \emph{oriented} if $\alpha$ is orientation-preserving or orientation-reversing. We denote the subsemigroups of $\mathcal{PT}_{n}$ consisting of all order-preserving transformations by $\mathcal{PO}_{n}$, consisting of all monotone transformations by $\mathcal{PM}_{n}$, consisting of all order-decreasing (extensive) transformations by $\mathcal{PD}_{n}$ ($\mathcal{PE}_{n}$), consisting of all orientation-preserving transformations by $\mathcal{POP}_{n}$, and consisting of all oriented transformations by $\mathcal{POR}_{n}$. We also let
\begin{eqnarray*}
	&&\mathcal{PC}_{n} = \mathcal{PO}_{n} \cap \mathcal{PD}_{n},\qquad\quad\,\,\,\, \mathcal{PMD}_{n} = \mathcal{PM}_{n} \cap \mathcal{PD}_{n},\\
	&&\mathcal{POPD}_{n} = \mathcal{POP}_{n} \cap \mathcal{PD}_{n},\quad \mathcal{PORD}_{n} = \mathcal{POR}_{n} \cap \mathcal{PD}_{n},\\
	&&\mathcal{IC}_{n} = \mathcal{PC}_{n} \cap \mathcal{I}_{n},\qquad\qquad\quad \mathcal{IOPD}_{n} = \mathcal{POPD}_{n} \cap \mathcal{I}_{n},\\
	&&\mathcal{IORD}_{n} = \mathcal{PORD}_{n} \cap \mathcal{I}_{n},
\end{eqnarray*}
and let $\mathcal{U}(n,r)=\{ \alpha\in \mathcal{U}_{n} : \lvert\im(\alpha)\rvert \leq r \}$ for every subsemigroup $\mathcal{U}_{n}$ of $\mathcal{PT}_{n}$ and $1\leq r\leq n$. It is shown in \cite[Lemma 1.1]{U2} that $\mathcal{D}_{n}= \mathcal{T}_{n} \cap \mathcal{PD}_{n}$ and $\mathcal{E}_{n} =\mathcal{T}_{n} \cap \mathcal{PE}_{n}$ are isomorphic semigroups. By a similar mapping to the one used in the proof of \cite[Lemma 1.1]{U2}, one can easily show that $\mathcal{PD}_{n}$ and  $\mathcal{PE}_{n}$, and so $\mathcal{PORD}_{n}$ and  $\mathcal{PORE}_{n}= \mathcal{POR}_{n} \cap \mathcal{PE}_{n}$ are isomorphic semigroups. Thus, it is enough to consider just $\mathcal{PORD}_{n}$.

We denote the smallest subsemigroup of a semigroup $S$ containing  a non-empty subset $A$ of $S$ by $\langle A\,\rangle$. If $\langle A\, \rangle =S$, then $A$ is called a \emph{generating set} of $S$. A generating set of $S$ with a minimum number of elements is called a \emph{minimal generating set}, the size of a minimal generating set of $S$ is called the \emph{rank} of $S$, and denoted by $\rank(S)$. An element $z$ of a semigroup $S$ is called \emph{undecomposable} in $S$ if, for all $x,y\in S$,\,  $z=xy$ implies $z=x$ or $z=y$. Then it is clear from the definition that every generating set of $S$ contains all undecomposable elements of $S$. The set of all idempotents in a non-empty subset $T$ of $S$ is denoted by $E(T)$, i.e. $E(T)=\{ e\in T: e^{2}=e\}$. It is clear that $\alpha\in \mathcal{PD}_{n}$ is an idempotent if and only if $\min(x\alpha^{-1}) =x$ for each $x\in \im(\alpha)$.  For other terms in semigroup theory, which are not explained here, we refer to \cite{GM,H}.

It is well-known that $\mathcal{PT}_{n}$, $\mathcal{T}_{n}$, and $\mathcal{I}_{n}$ have ranks $4$, $3$, and $3$, respectively. Gomes and Howie \cite{GH} showed that the ranks of the semigroup of all order-preserving full transformations $\mathcal{O}_n$ and order-preserving partial transformations $\mathcal{PO}_n$ are equal to $n+1$ and $2n$, respectively. In \cite{FGJ}, Fernandes, Gomes and Jesus showed that the rank of the monoid $\mathcal{POP}_n$ is equal to $3$. Dimitrova and Koppitz \cite{DK1} showed that the rank of $\mathcal{POE}_n = \mathcal{PO}_{n} \cap \mathcal{PE}_{n}$ is equal to $2n$. More recently, the rank of the semigroup $\mathcal{POPD}_{n}$ was determined by Ay\i k, Ay\i k, Bugay, and Da\u{g}deviren in \cite{AABD}. Ay\i k, Ay\i k, Dimitrova, and Koppitz determined the ranks of the semigroups $\mathcal{PMD}_{n}$, $\mathcal{IMD}_{n}$, and their ideals in \cite{AADK1}, as well as of the semigroup $\mathcal{ORD}_{n}$ and its ideals in \cite{AADK2}. Finding ranks of certain transformation semigroups is an important problem considered by many authors (see, for example, \cite{AAUH, D1, D2, F, G, LU2, LU3, LZL, ZP, ZHQ}). In the present paper, we will determine some minimal generating sets and ranks of $\mathcal{PORD}(n,r)$ and $\mathcal{IORD}(n,r)$ for $3\leq r\leq n$.

Moreover, we will characterize the maximal subsemigroups of $\mathcal{PORD}(n,r)$ as well as of $\mathcal{IORD}(n,r)$ for any $3\leq r\leq n$. By a maximal subsemigroup of a semigroup $S$, we mean a maximal element, under set inclusion, within the family of all proper subsemigroups of $S$. In 1966, Bayramov \cite{Bayramov} proved that $S$ is a maximal subsemigroup of $\mathcal{T}_{n}$ if and only if $S$ is either $\mathcal{T}(n,n-2) \cup \mathcal{S}_n$ or $\mathcal{T}(n,n-1) \cup G$, where $G$ is a maximal subgroup of the symmetric group $\mathcal{S}_n$. In 1968, Graham, Graham, and Rhodes \cite{GGR} described forms to which the maximal subsemigroups of a finite semigroup belong. Determining which of these forms appear in a given finite semigroup is difficult. Donoven, Mitchell, and Wilson \cite{DMW} presented an algorithm for computing the maximal subsemigroups of a finite semigroup. Since it is important to know the maximal subsemigroups for concrete semigroups, there are numerous papers about finding maximal subsemigroups of particular classes of transformation semigroups (see, for example, \cite{AADK1, AADK2, DFK, DK1, DK2, DM, EKMW, GM1, ZH, ZHQ}).

\section{Preliminaries}

The \emph{fix} and \emph{kernel} sets of $\alpha\in \mathcal{PT}_{n}$ are defined by $\fix(\alpha) =\{ x\in \dom(\alpha ): x\alpha=x\}$ and $\ker(\alpha) =\{ (x,y)\in \dom(\alpha)\times \dom(\alpha) : x\alpha= y\alpha\}$, respectively. We use the notations $1_{n}$ and $0_{n}$ to denote the identity and empty (zero) mappings on $X_{n}$, respectively. Let $Y$ be a non-empty subset of $X_{n}$. Then the restriction map of an element $\alpha$ in $\mathcal{PT}_{n}$ to $Y\cap \dom(\alpha)$ is denoted by $\alpha_{\mid_{Y}}$. In particular, the restriction map of $1_{n}$ to $Y$ is denoted by $1_{Y}$, and called the \emph{partial identity} on $Y$.

Let $U$, $V$ and $Y$ be any non-empty subsets of $X_{n}$. If $u<v$ for all $u\in U$ and $v\in V$, then we write $U<V$; and if $P=\{U_{1} ,\ldots, U_{r}\}$ ($1\leq r\leq n$) is a family of non-empty disjoint subsets of $Y$ such that $Y=\bigcup \limits_{i=1}^{r} U_{i}$, then $P$ is called a \emph{partition} of $Y$ with $r$ parts. For a non-zero $\alpha\in \mathcal{PT}_{n}$, notice that
$$\Ker(\alpha) =\{ y\alpha^{-1}: y\in \im(\alpha)\}$$
is a partition of $\dom(\alpha)$.  A partition $P=\{U_{1} ,\ldots, U_{r}\}$ of $Y$ is called an \emph{ordered} partition of $Y$ if $U_{i}< U_{i+1}$ for all $1\leq i\leq r-1$, and a \emph{convex} partition of $Y$ if for each $1\leq i\leq r$,\, $U_{i}$ is convex, i.e. for every $u,v\in U_{i}$ and $y\in Y$,\, $u\leq y\leq v$ implies $y\in U_{i}$. For a non-zero $\alpha\in \mathcal{PM}_{n}$, if $\im(\alpha) =\{x_{1}< \cdots< x_{r}\}$ ($1\leq r\leq n$), then it is clear that $\Ker(\alpha)$ is a convex and ordered partition of $\dom(\alpha)$. In addition, $\Ker(\alpha) =(x_{1} \alpha^{-1} <\cdots <x_{r} \alpha^{-1})$ if $\alpha\in \mathcal{PO}_{n}$, and otherwise, $\Ker(\alpha) =(x_{r} \alpha^{-1} <\cdots <x_{1} \alpha^{-1})$.

Let $\mathcal{PR}_{n}$ denote the subset of $\mathcal{PT}_{n}$ consisting of all orientation-reversing  partial transformations. It is shown in \cite[Lemma 1.1]{CH} that for $\alpha\in \mathcal{PT}_{n}$ with $\dom (\alpha) =\{a_{1}< \cdots < a_{t}\}$, the sequence $(a_{1}\alpha , \ldots ,a_{t}\alpha)$ is both cyclic and anti-cyclic if and only if $\lvert \im(\alpha) \rvert \leq2$. Thus, we have
\begin{eqnarray}\label{e1}
	\mathcal{POP}_{n}\cap\mathcal{PR}_{n} =\{\alpha\in \mathcal{POR}_{n} :\lvert \im(\alpha) \rvert\leq 2\}.
\end{eqnarray}
It is also shown in \cite[Proposition 2.3]{CH} that a non-constant oriented full transformation $\alpha$ is order-preserving (order-reversing) if and only if $1\alpha < n\alpha$ ($n\alpha < 1\alpha$). Similarly, it is shown that a non-constant $\alpha\in \mathcal{POR}_{n}$ is order-preserving (order-reversing) if and only if $(\min(\dom(\alpha))) \alpha < (\max(\dom(\alpha ))) \alpha$\, ($(\max(\dom(\alpha ))) \alpha <(\min(\dom(\alpha))) \alpha$).

We suppose that $n\geq 4$ since $\mathcal{PORD}_{n}= \mathcal{POPD}_{n} =\mathcal{PC}_{n}$ for $n=1,2$, and $\mathcal{PORD}_{3} =\mathcal{POPD}_{3}$, and that $3\leq r\leq n-1$ since $\mathcal{PORD}(n,1)= \mathcal{POPD}(n,1) =\mathcal{PC}(n,1)$,\, $\mathcal{PORD}(n,2)= \mathcal{POPD}(n,2)$,\, $\mathcal{PORD}(n,n-1)= \mathcal{PORD}_{n} \setminus \{ 1_{n}\}$ and $\mathcal{PORD}(n,n)= \mathcal{PORD}_{n}$ throughout this paper unless otherwise stated.

Next we let
\begin{eqnarray*}
	\mathcal{PMD}_{n}^{*}=\mathcal{PMD}_{n}\setminus \mathcal{PC}_{n},\qquad\,\,\,\,\,&&
	\mathcal{POPD}_{n}^{*}=\mathcal{POPD}_{n}\setminus \mathcal{PC}_{n},\\
	\mathcal{PORD}_{n}^{*}=\mathcal{PORD}_{n}\setminus \mathcal{POPD}_{n},&&
	\mathcal{PRD}_{n}^{*}=\mathcal{PORD}_{n}^{*}\setminus \mathcal{PMD}_{n},
\end{eqnarray*}
and for $1\leq r\leq s\leq n$, let
$$[r,s]=\{ r,r+1,\ldots, s\}.$$
In \cite[Section 2]{AABD}, the \emph{order-preserving degree} of $\alpha\in \mathcal{POPD}_{n}$ is defined by
$$\opd(\alpha) =\max\{ m : \alpha_{\mid_{X_{m}}} \in  \mathcal{PO}_{m} \}=\max\{ m : \alpha_{\mid_{X_{m}}} \in  \mathcal{PC}_{m}\},$$
and noticed that for $\alpha\in \mathcal{POPD}_{n}^{*}$, if $\opd(\alpha) =m$ and $\dom(\alpha) =\{x_{1} <\cdots< x_{t}\}$, then $2\leq t\leq n$,\, $2\leq m\leq n-1$, and there exists a unique $1\leq s\leq t-1$ such that $1\leq x_{s+1}\alpha \leq\cdots \leq x_{t}\alpha \leq x_{1}\alpha \leq\cdots \leq x_{s}\alpha \leq m$ and $x_{s+1}=m+1$. Thus, $\alpha_{\mid_{X_{m}}} \in \mathcal{PC}_{m}$, and if $k=x_{1}\alpha$, then $\alpha_{\mid_{[m+1,n]}}$ is an order-preserving partial transformation from  $[m+1,n]$ to $X_{k}$ such that $m+1\in \dom(\alpha_{\mid_{[m+1,n]}})$. Now we define the \emph{order-reserving degree} of $\alpha\in \mathcal{PORD}_{n}^{*}$ by
$$\ord(\alpha) =\max\{ m : \alpha_{\mid_{X_{m}}} \in \mathcal{PMD}_{n} \text{ and } \max(\im(\alpha))=(m+1)\alpha \}.$$
(This definition is a simple generalization of the definition given in \cite[Section 2]{AADK2}.) For $\alpha\in \mathcal{PRD}_{n}^{*}$, if $\ord(\alpha) =m$ and $\dom(\alpha) =\{x_{1} <\cdots< x_{t}\}$, then similarly, we notice that $3\leq t\leq n$,\, $2\leq m\leq n-1$, and there exists unique $1\leq s\leq t-1$ such that $x_{s+1} =m+1$ and $\im(\alpha) \subseteq [x_{s}\alpha, x_{s+1}\alpha]$. More precisely, $\alpha\in \mathcal{PORD}_{n}^{*}$ has the tabular form:
\begin{eqnarray}\label{e2}
	\alpha=\left(\begin{array}{ccc|ccc}
		x_{1}      &\cdots& x_{s}       & x_{s+1}       &\cdots& x_{t}\\
		x_{1}\alpha&\cdots& x_{s}\alpha & x_{s+1}\alpha &\cdots& x_{t}\alpha
	\end{array}\right)
\end{eqnarray}
with the property that
$$1\leq  x_{s}\alpha \leq\cdots\leq x_{1}\alpha \leq x_{t}\alpha \leq\cdots \leq x_{s+1}\alpha \leq x_{s+1}=m+1.$$
Thus, $\alpha_{\mid_{X_{m}}} \in \mathcal{PMD}_{m}$, and if $k=x_{1}\alpha$, then $\alpha_{\mid_{[m+1,n]}}$ is an order-reversing partial transformation from  $[m+1,n]$ to $[k,m+1]$ such that $m+1\in \dom(\alpha_{\mid_{[m+1,n]}})$.

Let $\alpha\in \mathcal{PRD}_{n}^{*}$ with $\ord(\alpha)=m$. As noticed in \cite[Equation (1)]{AADK1}, we have
\begin{eqnarray*}
	\lvert \im(\alpha_{\mid_{X_{m}}}) \rvert\leq \max(\im(\alpha_{\mid_{X_{m}}}))= (\min(\dom(\alpha_{\mid_{X_{m}}}))) \alpha\leq \min(\dom(\alpha_{\mid_{X_{m}}})),
\end{eqnarray*}
and so we conclude that
\begin{eqnarray}\label{e3}
	1\leq \lvert \im(\alpha_{\mid_{X_{m}}}) \rvert \leq \left \lfloor \frac{m+1}{2} \right \rfloor
\end{eqnarray}
where $\left \lfloor \frac{m+1}{2} \right \rfloor$ denotes the largest integer less than or equal to $\frac{m+1}{2}$. Notice that if
$\alpha=\left(\begin{matrix}
	3&5&7\\
	3&5&3
\end{matrix}\right)$
and
$\beta=\left(\begin{matrix}
	3&5&7\\
	3&5&4
\end{matrix}\right)$
in $\mathcal{PORD}_{7}$, then $\opd(\alpha)=6$ and $\ord(\beta)=4$, but $\ord(\alpha)$ and $\opd(\beta)$ are not defined since $\alpha \notin \mathcal{PORD}_{7}^{*}$ and $\beta \notin \mathcal{POPD}_{7}$.

Let $r_{n}$ denote the maximum size of image of a orientation-reversing and order-decreasing partial mapping on $X_{n}$, i.e.
$$r_{n}=\max\{ \lvert \im(\alpha) \rvert: \alpha\in \mathcal{PRD}_{n}^{*} \}.$$
Then we have the following result.

\begin{proposition}\label{p1}
	$r_{n}=n-\left \lfloor \frac{n}{3} \right \rfloor$.
\end{proposition}

\proof Let $k=\left \lfloor \frac{n}{3} \right \rfloor$, and so by the Division Algorithm, we have $n=3k+i$ for $0\leq i\leq 2$. If we consider the mapping
$$\gamma_{n} =\left( \begin{array}{cccc|cccc}
	k+i-1& k+i &\cdots&2k+2i-3&2k+i&2k+i+1&\cdots & n\\
	k+i-1&k+i-2&\cdots&   1   &2k+i&2k+i-1&\cdots &k+i
\end{array} \right),$$
then it is clear that $\gamma_{n} \in \mathcal{PRD}_{n}^{*}$ and $\lvert \im(\gamma_{n}) \rvert =2k+i= n-k$, and hence, $r_{n}\geq n- \left \lfloor \frac{n}{3} \right \rfloor$.

Let $\alpha\in \mathcal{PRD}_{n}^{*}$, and suppose that $\lvert \im(\alpha) \rvert =r\geq n-k$. If $\max(\im( \alpha))=t$, then there exists $0\leq j\leq k$ such that $\min(t\alpha ^{-1})= n-k+j= 2k+i+j$ since $\min(t\alpha ^{-1}) \geq t\geq r\geq n-k$, and so by (\ref{e2}), we obtain $m=\ord(\alpha) =2k+i+j-1$. Since $r=\lvert \im(\alpha) \rvert\leq  \lvert \im(\alpha _{\mid_{X_{m}}}) \rvert +\lvert \im(\alpha_{\mid_{[m+1,n]}}) \rvert$, it follows from (\ref{e3}) that
\begin{eqnarray*}
	r&\leq& \left \lfloor \frac{2k+i+j}{2} \right \rfloor +(k-j+1) =2k+i- \left( (j-1)- \left \lfloor \frac{j-i}{2} \right \rfloor \right)\\
	&=& 2k+i- \left \lfloor \frac{j-1+i}{2} \right \rfloor \leq 2k+i=n-k,
\end{eqnarray*}
and so $r=n-k$. Therefore, we conclude that $r_{n}=n-\left \lfloor \frac{n}{3} \right \rfloor$. \qed

For all $\delta, \lambda\in \mathcal{D}_{n}$, it is shown in \cite{LU2} that $\fix(\delta\lambda) =\fix(\delta) \cap\fix(\lambda)$. Also we give the following lemma without a proof since its proof is similar to the proof of \cite[Lemma 1.1]{LU2}.

\begin{lemma}\label{l2}
	$\fix(\delta\lambda) =\fix(\delta) \cap\fix(\lambda)$ for all $\delta, \lambda\in \mathcal{PD}_{n}$.\hfill $\square$
\end{lemma}

For every subsemigroup $\mathcal{U}_{n}$ of $\mathcal{PT}_{n}$ and $0\leq r\leq n$, let
$$E_{r}(\mathcal{U}_{n})=\{ \alpha \in E(\mathcal{U}_{n}) : \lvert\im(\alpha) \rvert =r \}.$$
For each $1\leq r\leq n$, it is determined in \cite[Proposition 1]{AABD} that $\lvert E_{r}(\mathcal{POPD}_{n})\rvert =\sum\limits_{s=r}^{n} \binom{n}{s} \binom{s}{r}=\binom{n}{r}2^{n-r}$. Now we have the following result.

\begin{proposition} \label{p3}
	For each $\alpha \in \mathcal{PORD}_{n}^{*}$, we have $\lvert \fix(\alpha) \rvert\leq 2$. Thus, $E(\mathcal{PORD}_{n}) =E(\mathcal{POPD}_{n})$, and so for each $1\leq r\leq n$, we have $\lvert E_{r}(\mathcal{PORD}_{n})\rvert =\binom{n}{r}2^{n-r}$.
\end{proposition}

\proof For any $\alpha\in \mathcal{PORD}_{n}^{*}$, let $\ord(\alpha)=m$ and $\dom(\alpha) =\{x_{1} <\cdots< x_{t}\}$. Then $1\leq  x_{s}\alpha \leq\cdots\leq x_{1}\alpha \leq x_{t}\alpha \leq\cdots \leq x_{s+1}\alpha \leq x_{s+1}=m+1$
for some $1\leq s\leq t-1$. Since $x_{1}\alpha \leq x_{1}$ and $\alpha$ is order-reversing on $\dom(\alpha) \cap X_{m}$, there exists at most one fixed point of $\alpha$ in $X_{m}$. Moreover, since $(m+1)\alpha \leq m+1$ and $\alpha$ is order-reversing on $\dom(\alpha) \cap [m+1,n]$, there exists at most one fixed point of $\alpha$ in $[m+1,n]$. Therefore, $\alpha$ has at most two fixed points in $X_{n}$.

For any $\alpha\in \mathcal{PORD}_{n}^{*}$, if $\alpha$ is an idempotent, then $\lvert \im(\alpha) \rvert =\lvert \fix(\alpha) \rvert \leq 2$, and so by (\ref{e1}), we obtain $\alpha\in E(\mathcal{POPD}_{n})$. Thus,  $E(\mathcal{PORD}_{n}) =E(\mathcal{POPD}_{n})$ and so $\lvert E_{r}(\mathcal{PORD}_{n}) \rvert= \binom{n}{r}2^{n-r}$ for each $1\leq r\leq n$.  \qed

\section{Minimal generating set and rank of $\mathcal{PORD}(n,r)$ and $\mathcal{IORD}(n,r)$}

For each $2\leq r\leq n-1$, it is proven in \cite[Proposition 5]{AABD} that each element of $E_{r}= E_{r}(\mathcal{POPD}_{n})$ is undecomposable in $\mathcal{POPD}(n,r)$. We also have the following proposition. Although the proof of \cite[Proposition 5]{AABD} is also valid for $\mathcal{PORD}(n,r)$, we want to give a short and easy proof to make it easier for the reader.

\begin{proposition}\label{p4}
	Each element of $E_{r}$ is undecomposable in $\mathcal{PORD}(n,r)$.
\end{proposition}

\proof For any $\xi\in E_{r}$, suppose that $\xi= \alpha \beta$ for some $\alpha, \beta\in \mathcal{PORD}(n,r)$. By Lemma \ref{l2}, since $\fix(\xi) =\im(\xi)$ is a subset of both $\fix(\alpha)$ and $\fix(\beta)$, we get $3\leq r\leq \lvert \fix(\alpha) \rvert,\, \lvert \fix(\beta)\rvert$, and so from Proposition \ref{p3}, both $\alpha$ and $\beta$ must be in $\mathcal{POPD}(n,r)$ which contradicts to the result given in \cite[Proposition 5]{AABD}. \qed

Next we consider some idempotents in $\mathcal{POPD}(n,r)$ whose image size are strictly less than $r$. For each $1\leq p\leq n-2$ and each $p+1\leq q\leq p+r-2\leq n-1$, we define
\begin{eqnarray}\label{e4}
	\xi_{p,q}^{r}&=&\left( \begin{matrix}
		p & p+1 & \cdots & q & q+1\\
		p & p+1 & \cdots & q & p
	\end{matrix}\right)  \in E(\mathcal{POPD}(n,r)),
\end{eqnarray}
and we	let
$$F_{r}=\{\, \xi_{p,q}^{r} :1\leq p\leq n-2\,\mbox{ and }\, p+1\leq q\leq p+r-2\leq n-1\, \}.$$
First notice that $\lvert \im(\xi_{p,q}^{r}) \rvert =q-p+1\leq r-1$, and so $E_{r}\cap F_{r}=\emptyset$. For each $3\leq r\leq n-1$, it is determined in \cite[Equation (3)]{AABD} that $\lvert F_{r}\rvert= \frac{(2n-r-1)(r-2)}{2}$. Moreover, it is proven in \cite[Proposition 6]{AABD} that each generating set of $\mathcal{POPD}(n,r)$ must contain a mapping such that the restriction to $[p,q+1]$ is $\xi_{p,q}^{r}$ for each $1\leq p\leq n-2$ and $p+1\leq q\leq p+r-2\leq n-1$. By using this fact, we have the following proposition.

\begin{proposition}\label{p5}
	Let $A$ be a generating set of $\mathcal{PORD}(n,r)$. Then, for each $1\leq p\leq n-2$ and each $p+1\leq q\leq p+r-2\leq n-1$, there exists an $\alpha \in A$ such that $\alpha_{\mid_{[p,q+1]}} =\xi_{p,q}^{r}$.
\end{proposition}

\proof Suppose that $A$ is a generating set of $\mathcal{PORD}(n,r)$. Let $1\leq p\leq n-2$ and $p+1\leq q\leq p+r-2\leq n-1$. Then, there exist $\alpha_{1} ,\ldots, \alpha_{t} \in A$ such that $\xi_{p,q}^{r} =\alpha_{t} \cdots \alpha_{1}$. Since $\lvert \fix(\xi_{p,q}^{r}) \rvert =q-p+1$, if $q\geq p+2$, then $\lvert \fix(\xi_{p,q}^{r}) \rvert =3$, and so it follows from Lemma \ref{l2} and Proposition \ref{p3} that $\alpha_{1},\ldots ,\alpha_{t} \in \mathcal{POPD}(n,r)$, and so from \cite[Proposition 6]{AABD}, we have $\alpha_{i \mid_{[p,q+1]}} =\xi_{p,q}^{r}$ for at least one $1\leq i\leq t$.

Suppose $q=p+1$. For each $1\leq i\leq t$, it follows from Lemma \ref{l2} that $p\alpha_{i} =p$ and $(p+1)\alpha_{i} =p+1$. Since $\Ker(\xi_{p,q}^{r}) =\{ \{p,p+2\}, \{p+1\}\}$ and each $\alpha_{i}$ is order-decreasing, it follows that either $(p+2)\alpha_{i} =p+2$ or $(p+2)\alpha_{i} =p$ for every $1\leq i\leq t$. Thus, we also have $\alpha_{i\mid_{[p,q+1]}}=\xi_{p,q}^{r}$ for at least one $1\leq i\leq t$. \qed

It is proven in \cite[Theorem 8]{AABD} that $E_{r}\cup F_{r}$ is a minimal generating set of $\mathcal{POPD}(n,r)$, and so $\rank(\mathcal{POPD}(n,r)) =\binom{n}{r} 2^{n-r}+ \frac{(2n-r-1)(r-2)}{2}$. In particular, $\rank(\mathcal{POPD}_{n}) =\frac{n^{2}+n+2}{2}$.

Now we focus on $\mathcal{PORD}_{n}^{*}$. Let $\gamma$ in $\mathcal{PORD}_{n}^{*}$ with $\fix(\gamma) =\{p,q\}$ where $1\leq p<q\leq n$. Then notice that $\dom(\gamma) \setminus \{p,q\} \neq \emptyset$ since $\lvert \dom(\gamma) \rvert\geq \lvert \im(\gamma) \rvert\geq 3$. Moreover, for each $t\in \dom(\gamma) \setminus \{p,q\}$, since $\gamma$ is order-decreasing, it is easy to see that $p<t$, and so
$$\dom(\gamma) \subseteq [p,n].$$
Also it is easy to see that $1\leq t\gamma \leq p\gamma =p$ for all $t\in \dom(\gamma) \cap [p,q-1]$, and that $1\leq p= p\gamma \leq t\gamma \leq q\gamma =q$ for all $t\in \dom(\gamma) \cap [q,n]$, and so
$$1\leq p\leq q-2\leq n-2,\quad \im(\gamma) \subseteq [1,q] \,\, \mbox{ and }\,\, \ord(\gamma) =q-1.$$
For every $1\leq p\leq q-2\leq n-2$, we let
$$k=\min\{p-1,q-p-1\}\,\,\mbox{ and }\,\, l=\min\{q-p-1,n-q\},$$
and define the mapping
\begin{eqnarray} \label{e5}
\gamma_{p,q}&=&\left(\begin{array}{cccc|cccc}
p&p+1&\cdots&p+k&q&q+1&\cdots&q+l\\
p&p-1&\cdots&p-k&q&q-1&\cdots&q-l
\end{array}\right) \in \mathcal{PORD}_{n}.
\end{eqnarray}
First notice that $p-k=1$ or $p+k=q-1$, and that $q-l=p+1$ or $q+l=n$. Thus, for every $\gamma\in \mathcal{PORD}_{n}^{*}$ with $\fix (\gamma)=\{p,q\}$, we have $\lvert \im(\gamma) \rvert \leq \lvert \im(\gamma_{p,q}) \rvert$. Also notice that
$\gamma_{1,n}=\left(\begin{array}{cc}
	1&n\\
	1&n
\end{array}\right) \in \mathcal{POPD}_{n}$. However, if we let
\begin{eqnarray*}
G_{n}=\{ \gamma_{p,q} : 1\leq p\leq q-2\leq n-2\}\setminus \{\gamma_{1,n}\},
\end{eqnarray*}
then is easy to see that $G_{n}$ is a subset of $\mathcal{PORD}_{n}^{*}$, and that
\begin{eqnarray*}
	\lvert G_{n}\rvert =-1+\sum\limits_{p=1}^{n-2}\sum\limits_{q=p+2}^{n} 1 = -1+\sum\limits_{p=1}^{n-2}(n-p-1)=\frac{n(n-3)}{2}.
\end{eqnarray*}
For example, $G_{4}=\{ \gamma_{1,3}, \gamma_{2,4}\}$ where $\gamma_{1,3}=\left(\begin{array}{ccc}
	1&3&4\\
	1&3&2
\end{array}\right)$ and $\gamma_{2,4}=\left(\begin{array}{ccc}
2&3&4\\
2&1&4
\end{array}\right)$.
As Proposition \ref{p5}, we also have the following proposition.

\begin{proposition}\label{p6}
	Let $A$ be a generating set of $\mathcal{PORD}(n,r)$. Then, for every $1\leq p\leq n-2$ and $p+2\leq q\leq n$, there exists an $\alpha \in A\cap \mathcal{PORD}^{*}(n,r)$ such that $\fix(\alpha) =\{p,q\}$.
\end{proposition}

\proof Let $A$ be a generating set of $\mathcal{PORD}(n,r)$. For every $1\leq p\leq q-2\leq n-2$, there exist $\alpha_{1} ,\ldots, \alpha_{t}\in A$ such that $\gamma_{p,q} =\alpha_{1} \cdots \alpha_{t}$. By Lemma \ref{l2}, we have $\{p,q\}\subseteq  \fix(\alpha_{i})$ for each $1\leq i\leq t$. Since the products of orientation-preserving mappings are also orientation-preserving, for at least one of $1\leq j\leq t$,\, $\alpha_{j}$ must be orientation-reversing. Since orientation-reversing mappings have at most two fix points, we conclude that $\fix(\alpha_{j}) =\{p,q\}$, as required. \qed

Next we determine which $\gamma_{p,q}$'s are undecomposable in $\mathcal{PORD}_{n}$ in the following lemma.

\begin{lemma}\label{l7}
	Let $1\leq p\leq q-2\leq n-2$. Then $\gamma_{p,q}$ is undecomposable in $\mathcal{PORD}_{n}$ if and only if $q-1,n\in \dom(\gamma_{p,q})$, i.e. $\dom(\gamma_{p,q}) =[p,n]$.
\end{lemma}

\proof $(\Rightarrow)$ Let $Y=\dom(\gamma_{p,q})$, and suppose that $q-1\notin Y$, i.e. $k=\min\{p-1, q-p-1\} =p-1< q-p-1$ as defined above. Then we observe that $p+k=2p-1\in Y$ and $(2p-1) \gamma_{p,q} =1$. If we let $Z_{1}= Y\cup \{q-1\}$ and consider  the partial identity $1_{Y}$ and the mapping $\beta_{1} :Z_{1} \rightarrow \im(\gamma_{p,q})$ defined by
$$x\beta_{1}=\left\{ \begin{array}{cl}
	x\gamma_{p,q} & \mbox{ if } x\in Y,\\
	1 & \mbox{ if } x=q-1,
\end{array}\right.$$
then $\gamma_{p,q}$ is different from both $1_{Y}$ and $\beta_{1}$, and $\gamma_{p,q} =1_{Y}\beta_{1}$. Thus, we conclude that $\gamma_{p,q}$ is not undecomposable in $\mathcal{PORD}_{n}$.

Suppose that $n\notin Y$, i.e. $l= \min\{q-p-1, n-q\} =q-p-1<n$ as defined above. Then we observe that $q+l= 2q-p-1\in Y$ and $(2q-p-1) \gamma_{p,q} =p+1$. If we let $Z_{2}= Y\cup \{n\}$ and consider the mapping $\beta_{2} :Z_{2} \rightarrow \im(\gamma_{p,q})$ defined by
$$x\beta_{2}=\left\{ \begin{array}{cl}
	x\gamma_{p,q} & \mbox{ if } x\in Y,\\
	p & \mbox{ if } x=n,
\end{array}\right.$$
then similarly, $\gamma_{p,q} =1_{Y}\beta_{2}$ is not undecomposable in $\mathcal{PORD}_{n}$.

$(\Leftarrow)$ Suppose that $q-1,n\in \dom(\gamma_{p,q})$, i.e. $\dom(\gamma_{p,q}) =[p,n]$. Assume that there are $\alpha, \beta\in \mathcal{PORD}_{n}$ such that $\gamma_{p,q} =\alpha \beta$. By Lemma \ref{l2}, $\{p,q\}$ is a subset of both $\fix(\alpha)$ and $\fix(\beta)$. First notice that $[p,n] \subseteq \dom(\alpha)$, and that $\alpha$ is an injective mapping on $[p,n]$ since $\gamma_{p,q}$ is injective.

Suppose that $\alpha$ is orientation-preserving, and so $\beta$ is orientation-reversing. Then $\alpha$ is order-preserving on $[p,q]$ otherwise $t\alpha<p<q$ for some $t\in (p,q)$ which contradicts with orientation-preservation of $\alpha$. Since $p=p\alpha \leq (p+1)\alpha \leq (p+1)$, it follows from the injectivity of $\alpha$ on $[p,n]$ that $(p+1)\alpha =(p+1)$. By continuing in this fashion, we obtain $(p+i)\alpha =(p+i)$ for each $1\leq i\leq q-p-1$. One can similarly show that $(q+j)\alpha =(q+j)$ for each $1\leq j\leq n-q$, and so we conclude that $\alpha_{\mid_{[p,n]}} =1_{p,n}$, and so $[p,n]\subseteq \dom(\beta)$ and $\beta_{\mid_{[p,n]}} =\gamma_{p,q}$. Since $\beta$ is orientation-reversing, we get $\dom(\beta) \subseteq [p,n]$. It follows that $\dom(\beta) =[p,n]= \dom(\gamma_{p,q})$, and so $\beta= \beta_{\mid_{[p,n]}} =\gamma_{p,q}$.

Suppose that $\alpha$ is orientation-reversing, and so $\beta$ is orientation-preserving. Then $\alpha$ is order-reversing on both $[p,q-1]$ and $[q,n]$. Since
$$p-1=((p+1)\alpha)\beta \leq (p+1)\alpha \leq p\alpha =p,$$
it follows from the injectivity of $\alpha$ on $[p,n]$ that $(p+1)\alpha =(p-1)$. By continuing in this fashion, we obtain $(p+i)\alpha =(p-i)$ for each $1\leq i\leq q-p-1$. One can similarly continue and prove that $\alpha= \gamma_{p,q}$, and so $\gamma_{p,q}$ is undecomposable in $\mathcal{PORD}_{n}$. \qed

If $n-\lfloor \frac{n}{3} \rfloor\leq r\leq n-1$, then by Proposition \ref{p1}, notice that $G_{n}$ of $\mathcal{PORD}(n,r)$. Now we state and prove one of the main results of this paper.

\begin{theorem}\label{t8}
	$E_{r}\cup F_{r}\cup G_{n}$ is a minimal generating set of $\mathcal{PORD}(n,r)$, and so $$\rank(\mathcal{PORD}(n,r)) =\binom{n}{r}2^{n-r}+ \frac{(2n-r-1)(r-2)}{2} +\frac{n(n-3)}{2}$$
	for each $n-\lfloor \frac{n}{3} \rfloor\leq r\leq n-1$.
\end{theorem}

\proof Let $n-\lfloor \frac{n}{3} \rfloor\leq r\leq n-1$. Since $\mathcal{POPD}(n,r) =\langle E_{r}\cup F_{r}\rangle$ and $\mathcal{PORD}_{n}^{*} =\mathcal{PORD}(n,r) \setminus \mathcal{POPD}(n,r)$, it is enough to show that $\mathcal{PORD}_{n}^{*} \subseteq \langle \mathcal{POPD}(n,r) \cup G_{n}\rangle$. Let $\alpha\in \mathcal{PORD}_{n}^{*}$ with $\ord(\alpha) =m$ and $\im(\alpha) =\{a_{1}, a_{2}, \ldots, a_{t}\}$ where $3\leq t\leq r$. Then we have the following sequence:
$$1\leq a_{s}< a_{s-1}< \cdots< a_{1}< a_{t}< a_{t-1}< \cdots< a_{s+1}\leq m+1.$$
for an unique $1\leq s\leq t-1$. If we take $A_{1}=\{ x\in X_{m}: x\alpha =a_{1}\}$,\, $A_{i}=a_{i} \alpha^{-1}$ for each $2\leq i\leq t$ and $A_{t+1}=\{ x\in [m+1,n]: x\alpha =a_{1}\}$ ($A_{t+1}$ may be the empty set), then $\alpha$ can be written in the following form:
$$\alpha =\left(\begin{array}{cccc|cccc}
	A_{1} & A_{2} & \cdots & A_{s} & A_{s+1} & \cdots & A_{t} & A_{t+1}\\
	a_{1} & a_{2} & \cdots & a_{s} & a_{s+1} & \cdots & a_{t} &  a_{1}
\end{array}\right)$$
with the property that either
\begin{itemize}
	\item $(A_{1},\ldots, A_{s}, A_{s+1},\ldots, A_{t})$, whenever $A_{t+1}=\emptyset$, or
	\item $(A_{1}, \ldots, A_{s}, A_{s+1}, \ldots, A_{t}, A_{t+1})$, whenever $A_{t+1}\neq \emptyset$,
\end{itemize}
is a convex and ordered partition of $\dom(\alpha)$. In addition, since $\alpha$ is order-decreasing and $\min(A_{s+1})=m+1$, if we take $a_{1}=p$ to simplify notation, then we observe that
\begin{itemize}
	\item [$(i)$]   $a_{i}\leq p-i+1$ for every $1\leq i\leq s$,
	\item [$(ii)$] $p+i-1\leq \min(A_{i})$ for every $1\leq i\leq s$,
	\item [$(iii)$]  $a_{s+j}\leq m+2-j$  for every $1\leq j\leq t-s$, and
	\item [$(iv)$]  $m+j\leq \min(A_{s+j})$ for every $1\leq j\leq t-s$.
\end{itemize}
After all these observations, we define the mappings
$$\beta=\left(\begin{array}{cccccccc|c}
	A_{1} & A_{2} &\cdots & A_{s} & A_{s+1} & A_{s+2} &\cdots & A_{t} & A_{t+1}\\
	p   &  p+1  &\cdots & p+s-1 &   m+1   &   m+2   &\cdots & m+t-s &   p
\end{array}\right)$$
and
$$\delta=\left(\begin{array}{cccccccc}
	p-s+1 &\cdots &  p-1  &   p   &  m-t+s+2  &\cdots &    m    &  m+1\\
	a_{s} &\cdots & a_{2} & a_{1} &   a_{t}   &\cdots & a_{s+2} & a_{s+1}
\end{array}\right).$$
By $(ii)$, since $p+s-1\leq \min(A_{s})\leq m$, it follows from $(ii)$ and $(iv)$ that $\beta\in \mathcal{POPD}(n,r)$. By $(iii)$, since $p=a_{1}< a_{t}\leq m-t+s+2$, i.e. $p< m-t+s+2$, it follows from $(i)$ and $(iii)$ that $\delta\in \mathcal{PC}(n,r) \subseteq \mathcal{POPD}(n,r)$. Also we define
\begin{equation}\label{e6}
	\gamma =\left(\begin{array}{cccc|cccc}
		p & p+1 &\cdots & p+s-1 & m+1 & m+2 &\cdots & m+t-s\\
		p & p-1 &\cdots & p-s+1 & m+1 &  m  &\cdots & m-t+s+2
	\end{array}\right).
\end{equation}
Since $\lvert \im(\gamma) \rvert=t\geq 3$, it is clear that $\gamma\in \mathcal{PORD}_{n}^{*}$ and $\gamma$ is an injective mapping with $\fix(\gamma) =\{p,m+1\}$, and that $\alpha =\beta \gamma \delta$. If $\gamma \in G_{n}$, then the proof is completed. Let $Y=\dom(\gamma)$. Since $p+s-1\leq \min(A_{s})\leq m$ and $1\leq a_{s}\leq p-s+1$, we have $s-1\leq \min\{p-1,m-p\}$, and since $m+t-s\leq \min(A_{t}) \leq n$ and $p=a_{1}< a_{t}\leq m-t+s+2$, we have $t-s-1\leq \min\{n-(m+1), (m+1)-p-1\}$. Thus, it follows from the definition of $\gamma_{p,m+1}$ which is given in (\ref{e5}) that $Y\subseteq \dom(\gamma_{p,m+1})$, and that $\gamma= 1_{Y}\gamma_{p,m+1}$ which proves that $\gamma\in \langle \mathcal{POPD}(n,r) \cup G_{n}\rangle =\langle E_{r}\cup F_{r} \cup G_{n}\rangle$, and so $\mathcal{PORD}(n,r) =\langle E_{r}\cup F_{r} \cup G_{n}\rangle$.

Finally, it follows from Propositions \ref{p4}, \ref{p5} and \ref{p6} that $E_{r}\cup F_{r} \cup G_{n}$ is a minimal generating set. Moreover, since $E_{r}$, $F_{r}$ and $G_{n}$ are disjoint subsets, we have
$\rank(\mathcal{PORD}(n,r)) =\binom{n}{r}2^{n-r}+ \frac{(2n-r-1)(r-2)}{2} +\frac{n(n-3)}{2}$,
as claimed. \qed

Since $\rank(\mathcal{POPD}_{n}) =\frac{n^{2}+n+2}{2}$, we also have the following immediate corollary.

\begin{corollary}\label{c9}
	$\rank(\mathcal{PORD}_{n}) =n^{2}-n+1$. \hfill $\square$
\end{corollary}

If $p=n-2\lfloor \frac{n}{3} \rfloor-1$ and $q=n-\lfloor \frac{n}{3} \rfloor$, then, since $\lvert \im(\gamma_{p,q}) \rvert=n-\lfloor \frac{n}{3} \rfloor$, we observe that $G_{n}$ is not a subset of $\mathcal{PORD}(n,r)$, whenever $1\leq r<n-\lfloor \frac{n}{3} \rfloor$.

Let $3\leq r<n-\lfloor \frac{n}{3} \rfloor$ and let $1\leq p\leq q-2\leq n-2$. If $\alpha\in \mathcal{PORD}(n,r)$ is injective with $\fix(\alpha) =\{p,q\}$ and if we take $s=\lvert \dom(\alpha) \cap [p,q-1] \rvert$ and $u=\lvert \dom(\alpha) \cap [q,n] \rvert$, then we notice that $s\leq \min\{p,q-p,r-1\}$ and $u\leq \min\{r-s,q-p,n-q+1\}$. Now, for each $1\leq s\leq \min\{p,q-p,r-1\}$, we let $u=\min\{r-s,q-p,n-q+1\}$ and define the mapping
\begin{equation}\label{e7}
	\gamma_{p,q}^{r,s} =\left(\begin{array}{cccc|cccc}
		p & p+1 &\cdots & p+s-1 & q & q+1 &\cdots & q+u-1\\
		p & p-1 &\cdots & p-s+1 & q & q-1 &\cdots & q-u+1
	\end{array}\right) .
\end{equation}
Since $\lvert \im(\gamma_{p,q}^{r,s}) \rvert =s+u\leq s+(r-s)=r$, we notice that $\gamma_{p,q}^{r,s} \in \mathcal{PORD}(n,r)$. For each $1\leq p\leq n-2$, we also notice that $\gamma_{p,n}^{r,1}=\left(\begin{array}{c|c}
	p&n\\
	p&n
\end{array}\right) \in \mathcal{POPD}(n,r)$. However, if we let
$$H_{n}^{r}=\{ \gamma_{p,q}^{r,s} : 1\leq p\leq q-2\leq n-2,\,\,\, 1\leq s\leq \min\{p,q-p,r-1\} \mbox{ and } (q,s)\neq (n,1)\},$$
then it is easy to see that $H_{n}^{r}$ is a subset of $\mathcal{PORD}_{n}^{*}$ since $u\geq 2$, whenever $s=1$. Furthermore, since $p+1= q-(q-p) +1\leq q-u+1$, by (\ref{e7}), it is also easy to see that every element of $H_{n}^{r}$ is injective, that is $H_{n}^{r}$ is a subset of $\mathcal{IORD}(n,r)$. With these notations, we have the following proposition:

\begin{proposition}\label{p10}
	Let $3\leq r<n-\lfloor \frac{n}{3}\rfloor$ and $1\leq p\leq q-2\leq n-2$, and let $A$ be a generating set of $\mathcal{PORD}(n,r)$. Then, for each $1\leq s\leq \min\{p,q-p,r-1\}$, if $(q,s)\neq (n,1)$, then there exists an $\alpha\in A$ such that $\alpha_{\mid_{Y}} =\gamma_{p,q}^{r,s}$, where $Y=\dom(\gamma_{p,q}^{r,s})$.
\end{proposition}

\proof Let $1\leq s\leq \min\{p,q-p,r-1\}$ with $(q,s)\neq (n,1)$. Then, there exist $\alpha_{1},\ldots ,\alpha_{k} \in A$ such that $\gamma_{p,q}^{r,s} =\alpha_{1} \cdots \alpha_{k}$. By Lemma \ref{l2}, we have $\{p,q\}\subseteq  \fix(\alpha_{i})$ for each $1\leq i\leq k$. Moreover, note that $(q,s)\neq (n,1)$ justifies that $\gamma_{p,q}^{r,s}$ is orientation-reserving but not orientation-preserving. So, since the products of orientation-preserving mappings are also orientation-preserving, there exists $1\leq j\leq k$ such that $\alpha_{j}$ is orientation-reversing, and so by Proposition \ref{p3}, $\fix(\alpha_{j}) =\{p,q\}$. Without loss of generality, we suppose that $\beta=\alpha_{1} \cdots \alpha_{j-1}$ is orientation-preserving. (We take $\beta=1_{n}$ if $j=1$.) Since $\gamma_{p,q}^{r,s}$ is order-decreasing and injective, it follows that the restriction of $\beta$ to $Y=\dom(\gamma_{p,q}^{r,s})$ is order-decreasing, orientation-preserving and injective, and so we immediately have $\beta_{\mid_Y} =1_{Y}$. Thus, we also have the equality $\gamma_{p,q}^{r,s} =\alpha_{j} \cdots \alpha_{k}$. Since $\dom(\alpha_{j}) \subseteq [p,n]$ and $\fix(\alpha_{j}) =\{p,q\}$, it is clear that $\alpha_{j}$ is order-reversing on both $[p,p+s-1]\subseteq \dom(\alpha) \cap [p,q-1]$ and $[q,q+u-1]\subseteq \dom(\alpha) \cap [q,n]$. Since $\alpha_{j}$ is also injective on $Y$, it follows that $p-1\leq (p+1)\alpha_{j} <p\alpha_{j}=p$, i.e. $(p+1)\alpha_{j}=p-1$, whenever $s\geq 2$. If we continue in this fashion, we obtain the equality $\alpha_{j\mid_{Y}} =\gamma_{p,q}^{r,s}$, as claimed. \qed

Now we have the following similar theorem to Theorem \ref{t8}.

\begin{theorem}\label{t11}
	Let $3\leq r<n-\lfloor \frac{n}{3} \rfloor$. Then $E_{r}\cup F_{r}\cup H_{n}^{r}$ is a minimal generating set of $\mathcal{PORD}(n,r)$
	for each $3 \leq r < n-\lfloor\frac{n}{3}\rfloor$.
\end{theorem}

\proof With the notations and assumptions given in the proof of Theorem \ref{t8}, it is enough to show that $\gamma\in \langle \mathcal{POPD}(n,r) \cup H_{n}^{r}\rangle$ since $\alpha =\beta \gamma \delta$. As defined in (\ref{e6}) and (\ref{e7}), it is clear that $1\leq s\leq \min\{p,m+1-p,r-1\}$ and that $t-s\leq u=\min\{r-s, m+1-p, n-m\}$ since $3\leq t\leq r$, and so $Y=\dom(\gamma) \subseteq \dom(\gamma_{p,m+1}^{r,s})$. Now it is clear from their definitions that $\gamma= 1_{Y} \gamma_{p,m+1} ^{r,s}$, and so $\mathcal{PORD}(n,r) =\langle E_{r}\cup F_{r}\cup H_{n}^{r}\rangle$. Therefore, by Propositions \ref{p4}, \ref{p5} and \ref{p10}, $E_{r}\cup F_{r} \cup H_{n}^{r}$ is a minimal generating set. \qed

Recall that $E(\mathcal{I}_{n})=\{ 1_{Y} : \emptyset \neq Y\subseteq X_{n}\} \cup \{0_{n}\}$ is the set of all idempotents of $\mathcal{I}_{n}$. For $1\leq r\leq n-1$ and $2\leq a\leq n$, let $Y$ be any subset of $X_{n}\setminus \{ a-1,a\}$ with size $r-1$ and let $\delta^{a}_{Y} : Y\cup\{a\}\rightarrow X_{n}$ be defined by
$$x\delta^{a}_{Y}=\left\{ \begin{array}{lcl}
	x & \mbox{ if } & x\in Y\\
	a-1 & \mbox{ if } & x=a
\end{array}\right. $$
for all $x\in Y\cup\{a\}$. Moreover, for $1\leq r\leq n-1$, let
\begin{eqnarray*}
	EI_{r} &=& \{ 1_{Y} : Y \mbox{ is a subset of } X_{n} \mbox{ with size } r\} \mbox{ and let} \\
	FI_{r} &=& \{ \delta^{a}_{Y} : 2\leq a\leq n\, \mbox{ and }\, Y \subseteq X_{n}\setminus \{ a-1,a\} \mbox{ with } \lvert Y\rvert =r-1\}.
\end{eqnarray*}
For $1\leq r\leq n-1$, it is proven in \cite[Theorem 4]{AAD} that $EI_{r}\cup FI_{r}$ is the unique minimal generating set of $\mathcal{IC}(n,r)$, and that $\rank(\mathcal{IC}(n,r))=\binom{n}{r}+r \binom{n-1}{r}$. In particular, it is concluded in \cite[Corollary 5]{AAD} that $\rank(\mathcal{IC}_{n})=2n$.

For $2\leq k\leq n-1$, let $Z$ be any subset of $X_{n}\setminus\{1,n\}$ with size $k-1$. If $\max(Z)=a$ and $\min(Z)=b$, then let $\zeta_{Z} : Z\cup\{a+1\}\rightarrow X_{n}$ be defined by
$$x\zeta_{Z}=\left\{ \begin{array}{lcl}
	x & \mbox{ if } & x\in Z\\
	b-1 & \mbox{ if } & x=a+1
\end{array}\right. $$
for all $x\in Z\cup\{a+1\}$. Moreover, for $2\leq k\leq n-1$, let
\begin{eqnarray*}
	&&GI_{k}=\{ \zeta_{Z} :Z \mbox{ is a subset of } X_{n}\setminus \{ 1,n\} \mbox{ with size } k-1\}\, \mbox{ and let}\\
	&&GI^{c}_{k}=\{ \zeta_{Z} :Z \mbox{ is a convex subset of } X_{n}\setminus \{ 1,n\} \mbox{ with size }k-1 \},
\end{eqnarray*}
which is a subset of $GI_{k}$. For $1\leq r\leq n-1$, it is proven in \cite[Theorem 8]{AAD} that $EI_{r}\cup FI_{r}\cup GI_{r}\cup \left( \bigcup\limits_{k=2}^{r-1}GI^{c}_{k}\right)$ is the unique minimal generating set of $\mathcal{IOPD}(n,r)$, and that $\rank(\mathcal{IOPD}(n,r))= \binom{n}{r}+ n\binom{n-2}{r-1} +(r-2)n- \frac{r^{2}-r-2}{2}$. In particular, it is concluded in \cite[Corollary 9]{AAD} that $\rank(\mathcal{IOPD}_{n}) =\frac{n^{2}+n+2}{2}$.

Now we state the following lemma without a proof since the proofs of \cite[Lemmas 3 and 7]{AAD} are also valid in $\mathcal{IORD}(n,r)$.

\begin{lemma}\label{il2}
	\begin{itemize}
		\item [$(i)$] For $1\leq r\leq n-1$, each element of $EI_{r}\cup FI_{r}$ is undecomposable in $\mathcal{IORD}(n,r)$.
		\item [$(ii)$] For $2\leq r\leq n-1$, each element of $GI_{r}\cup \left( \bigcup\limits_{k=2}^{r-1}GI^{c}_{k}\right)$ is  undecomposable in $\mathcal{IORD}(n,r)$. \hfill $\square$
	\end{itemize}
\end{lemma}

Next we let
\begin{eqnarray*}
	\mathcal{IORD}_{n}^{*}=\mathcal{IORD}_{n}\setminus \mathcal{IOPD}_{n}.
\end{eqnarray*}

Since the proof of Proposition \ref{p6} is also valid in $\mathcal{IORD}(n,r)$, we just state the following proposition.

\begin{proposition}\label{ip3}
	Let $A$ be a generating set of $\mathcal{IORD}(n,r)$. Then, for every $1\leq p\leq n-2$ and $p+2\leq q\leq n$, there exists an $\alpha \in A\cap \mathcal{IORD}^{*}(n,r)$ such that $\fix(\alpha) =\{p,q\}$. \hfill $\square$
\end{proposition}

Moreover, we have the following lemma.

\begin{lemma}\label{il4}
	For $1\leq p\leq n-2$ and $p+2\leq q\leq n$,\, $\gamma_{p,q}$ is undecomposable in $\mathcal{IORD}_{n}$.
\end{lemma}

\proof Let $1\leq p\leq q-2\leq n-2$ and let $Y=\dom(\gamma_{p,q})$. If $q-1,n\in Y$, then it follows from  Lemma \ref{l7} that $\gamma_{p,q}$ is undecomposable in $\mathcal{PORD}_{n}$, and so in $\mathcal{IORD}_{n}$. Thus, it is enough to consider the cases: $q-1\notin Y$ or $n\notin Y$.

Suppose that $q-1\notin Y$, i.e. $k=\min\{p-1, q-p-1\} =p-1< q-p-1$ as defined in (\ref{e5}). Thus, $p+k=2p-1\in Y$, i.e. $Y\cap [p,q-1] =[p,2p-1]$ and $(2p-1) \gamma_{p,q} =1$. Then assume that there are $\alpha, \beta\in \mathcal{IORD}_{n}$ such that $\gamma_{p,q} =\alpha \beta$. By Lemma \ref{l2}, we have $\fix(\gamma_{p,q}) =\{p,q\} \subseteq \fix(\alpha) \cap\fix(\beta)$. Moreover, notice that $Y\subseteq \dom(\alpha)$, and that $\alpha$ is injective on $Y$ since $\gamma_{p,q}$ is injective. Suppose $\alpha$ is orientation-preserving, then $\beta$ is orientation-reversing. In fact, as shown in the proof of Lemma \ref{l7}, it is similarly shown that $\alpha_{\mid_{Y}} =1_{Y}$, and so $Y\subseteq \dom(\beta)$ and $1_{Y}\beta =\gamma_{p,q}$. Since $\beta\in \mathcal{IORD}_{n}^{*}$ with $\fix (\beta)=\{p,q\}$, it follows from the definition of $\gamma_{p,q}$, given in (\ref{e5}), that $\lvert \dom(\beta) \rvert =\lvert \im(\beta) \rvert \leq \lvert \im(\gamma_{p,q}) \rvert =\lvert Y\rvert$, and so $\dom(\beta)=Y$, Thus, we have $\gamma_{p,q} =\beta$.

Suppose that $\alpha$ is orientation-reversing. Since $\fix(\alpha) =\{p,q\}$, it similarly follows that $\dom(\alpha)=Y$. Since the mappings are injective, as shown in the proof of Lemma \ref{l7}, we have $\gamma_{p,q} =\alpha$. Thus, $\gamma_{p,q}$ is undecomposable in  $\mathcal{IORD}_{n}$.

Since the proof for $n\notin Y$ is dual to the proof for $q-1\notin Y$, the proof is completed.\qed

With above notations, we are now able to give the following result:

\begin{theorem}\label{it5}
	$EI_{r}\cup FI_{r}\cup GI_{r}\cup \left( \bigcup\limits_{k=2}^{r-1}GI^{c}_{k}\right)\cup G_{n}$ is the unique minimal generating set of $\mathcal{IORD}(n,r)$, and so $$\rank(\mathcal{IORD}(n,r)) =\binom{n}{r}+ n\binom{n-2}{r-1}+(r-2)n-\frac{r^{2}-r-2}{2}+\frac{n(n-3)}{2}$$
	for each $n-\lfloor \frac{n}{3} \rfloor\leq r\leq n-1$. In particular, $\rank(\mathcal{IORD}_{n}) =n^{2}-n+1$.
\end{theorem}

\proof Let $n-\lfloor \frac{n}{3} \rfloor\leq r\leq n-1$. Since $\mathcal{IOPD}(n,r) =\langle EI_{r}\cup FI_{r}\cup GI_{r}\cup \left( \bigcup\limits_{k=2}^{r-1}GI^{c}_{k}\right)\rangle$, it is enough to show that $\mathcal{IORD}_{n}^{*} \subseteq \langle \mathcal{IOPD}(n,r) \cup G_{n}\rangle$. Let $\alpha\in \mathcal{IORD}_{n}^{*}$ with $\ord(\alpha) =m$ and $\im(\alpha) =\{a_{1}, a_{2}, \ldots, a_{t}\}$ where $3\leq t\leq r$. Then we have the following sequence:
$$1\leq a_{s}< a_{s-1}< \cdots< a_{1}< a_{t}< a_{t-1}< \cdots< a_{s+1}\leq m+1.$$
for an unique $1\leq s\leq t-1$. If we take $b_{i}=a_{i} \alpha^{-1}$ for each $1\leq i\leq t$, then $\alpha$ can be written in the following tabular form:
$$\alpha =\left(\begin{array}{cccc|ccc}
	b_{1} & b_{2} & \cdots & b_{s} & b_{s+1} & \cdots & b_{t} \\
	a_{1} & a_{2} & \cdots & a_{s} & a_{s+1} & \cdots & a_{t}
\end{array}\right) $$
with the property that $1\leq b_{1}< \cdots <b_{s} < b_{s+1}=m+1 < b_{s+2}< \cdots <b_{t}$. In addition, if we take $a_{1}=p$, then we observe that
\begin{itemize}
	\item [$(i)$]   $a_{i}\leq p-i+1$ for every $1\leq i\leq s$,
	\item [$(ii)$] $p+i-1\leq b_{i}$ for every $1\leq i\leq s$,
	\item [$(iii)$]  $a_{s+j}\leq m+2-j$  for every $1\leq j\leq t-s$, and
	\item [$(iv)$]  $m+j\leq b_{s+j}$ for every $1\leq j\leq t-s$.
\end{itemize}
After all these observations, if we define the mapping
$$\beta=\left(\begin{array}{cccccccc}
	b_{1} & b_{2} &\cdots & b_{s} & b_{s+1} & b_{s+2} &\cdots & b_{t} \\
	p   &  p+1  &\cdots & p+s-1 &   m+1   &   m+2   &\cdots & m+t-s
\end{array}\right)$$
and consider $\delta$ as defined in the proof of Theorem \ref{t8}, then it similarly follows from $(ii)$ and $(iv)$ that $\beta\in \mathcal{IC}(n,r)$, and it similarly follows from $(i)$ and $(iii)$ that $\delta\in \mathcal{IC}(n,r)$. Moreover, if we consider $\gamma$ as defined in (\ref{e6}), then it is clear that $\gamma\in \mathcal{IORD}_{n}^{*}$,\, $\fix(\gamma) =\{p,m+1\}$, and $\alpha =\beta \gamma \delta$. Thus, it similarly follows from the definition of $\gamma_{p,m+1}$, which is given in (\ref{e5}), Lemmas \ref{il2} and \ref{il4} that $EI_{r}\cup FI_{r}\cup GI_{r}\cup \left( \bigcup\limits_{k=2}^{r-1}GI^{c}_{k}\right) \cup G_{n}$ is the unique minimal generating set of $\mathcal{IORD}(n,r)$, and so by \cite[Theorem 8]{AAD}, we have
$$\rank(\mathcal{IORD}(n,r))= \binom{n}{r}+ n\binom{n-2}{r-1}+(r-2)n-\frac{r^{2}-r-2}{2}+\frac{n(n-3)}{2}.$$
In particular, $\rank(\mathcal{IORD}_{n})= \frac{n^{2}+n+2}{2} +\frac{n(n-3)}{2} =n^{2}-n+1$ since $\mathcal{IORD}_{n}= \mathcal{IORD}(n,n-1) \cup \{1_{n}\}$ and $1_{n}$ is undecomposable in $\mathcal{IORD}_{n}$. \qed

Let $3\leq r<n-\lfloor \frac{n}{3} \rfloor$ and let $1\leq p\leq q-2\leq n-2$. For each $1\leq s\leq \min \{p,q-p, r-1\}$, consider $\gamma_{p,q}^{r,s}$, as defined in (\ref{e7}), and $H_{n}^{r}$. Since the proofs of Proposition \ref{p10} and Theorem \ref{t11} are also valid for the following theorem, we omit its proof.

\begin{theorem}\label{it6}
	Let $3\leq r<n-\lfloor \frac{n}{3} \rfloor$. Then $EI_{r}\cup FI_{r}\cup GI_{r}\cup \left( \bigcup\limits_{k=2}^{r-1}GI^{c}_{k}\right)\cup H_{n}^{r}$ is a minimal generating set of $\mathcal{IORD}(n,r)$. \hfill $\square$
\end{theorem}

Notice that $\gamma_{5,7}^{5,1}=\left(\begin{array}{c|cc}
	5 & 7 & 8\\
	5 & 7 & 6
\end{array}\right) \in \mathcal{IORD}(9,5)$ is not undecomposable in $\mathcal{IORD}(9,5)$ since $1_{Y} \left(\begin{array}{cc|cc}
	5 & 6 & 7 & 8\\
	5 & 4 & 7 & 6
\end{array}\right) =\gamma_{5,7}^{5,1}$, where $Y=\{5,7,8\}$.
However, we have the following lemma.

\begin{lemma}\label{il7}
	Let $3\leq r<n-\lfloor \frac{n}{3} \rfloor$,\, $1\leq p\leq n-2$ and $p+2\leq q\leq n$. For each $1\leq s\leq \min\{p,q-p,r-1\}$, if $\min\{r-s,q-p,n-q+1\} =r-s$, then  $\gamma_{p,q}^{r,s}$ is undecomposable in $\mathcal{IORD}(n,r)$.
\end{lemma}

\proof We have $1\leq p\leq q-2\leq n-2$. Let $Y=\dom(\gamma_{p,q}^{r,s})$ and $Z=\im(\gamma_{p,q}^{r,s})$. For each $1\leq s\leq \min\{p,q-p,r-1\}$, if $\min\{r-s,q-p,n-q+1\} =r-s$, then by (\ref{e7}), we first notice that $(q,s)\neq (n,1)$, and that
$$\gamma_{p,q}^{r,s} =\left(\begin{array}{cccc|cccc}
	p & p+1 &\cdots & p+s-1 & q & q+1 &\cdots & q+r-s-1\\
	p & p-1 &\cdots & p-s+1 & q & q-1 &\cdots & q-r+s+1
\end{array}\right),$$
and so $\lvert \dom(\gamma_{p,q}^{r,s})\rvert =r$. Assume that there are $\alpha, \beta\in \mathcal{IORD}(n,r)$ such that $\gamma_{p,q}^{r,s} =\alpha \beta$. Then it is clear that $\dom(\alpha) =Y$ and $\im(\beta) =Z$. Then, after all these experiences, one can easily show that either $\alpha =1_{Y}$ and $\beta= \gamma_{p,q}^{r,s}$, or $\alpha= \gamma_{p,q}^{r,s}$ and $\beta =1_{Z}$.\qed

\section{Maximal subsemigroups of $\mathcal{PORD}(n,r)$ and\\ $\mathcal{IORD}(n,r)$}

In this section, we will characterize the maximal subsemigroups of $\mathcal{PORD}(n,r)$ as well as of  $\mathcal{IORD}(n,r)$ for $3 \leq r \leq n$.

For $1\leq p\leq n-2$ and $p+1\leq q\leq p+r-2\leq n-1$,  let
$$F^{r}_{p,q}=\{\, \alpha \in \mathcal{PORD}(n,r) : \alpha_{\mid_{[p,q+1]}} =\xi_{p,q}^{r} \}$$
and for $1\leq p\leq q-2 \leq n-2$ with $(p,q) \neq (1,n)$, let
$$G_{p,q} =\{ \alpha \in \mathcal{PORD}^*(n,r) : \alpha_{\mid_{\dom(\gamma_{p,q})}} = \gamma_{p,q} \}.$$

If $3 \leq r < n-\lfloor\frac{n}{3}\rfloor$, then let
$$H_{p,q}^{r,s} = \{\alpha \in \mathcal{PORD}^*(n,r) : \alpha_{\mid_{ \dom(\gamma_{p,q}^{r,s})}} = \gamma_{p,q}^{r,s} \}$$
for $1\leq p\leq q-2 \leq n-2$ and $1 \leq s \leq \min\{p, q-p, r-1\}$ with $(q,s) \neq (n,1)$.

\begin{lemma}\label{l12}
	Let $n-\lfloor\frac{n}{3}\rfloor \leq r \leq n-1$,\, $1 \leq p \leq n-2$, and $p+1 \leq q \leq p+r-2\leq n-1$. Then $\mathcal{PORD}(n,r)\setminus F^{r}_{p,q}$ is a maximal subsemigroup of $\mathcal{PORD}(n,r)$.
\end{lemma}

\proof Assume that there are $\alpha, \beta\in \mathcal{PORD}(n,r) \setminus F^{r}_{p,q}$ such that $\alpha \beta\in F^{r}_{p,q}$. Then $\alpha \beta_{\mid_{[p,q+1]}} =\xi_{p,q}^{r}$, i.e. $1_{[p,q+1]} \alpha \beta= \xi_{p,q}^{r}$, and so by Lemma \ref{l2}, we have $x\alpha=x$ and $x\beta=x$ for all $x\in[p,q]$. Therefore, we obtain $\alpha_{\mid_{[p,q]}} = \beta_{\mid_{[p,q]}} = \xi_{p,q\mid_{[p,q]}}$. From $p=(q+1) \xi_{p,q}^{r} =(q+1) 1_{[p,q+1]} \alpha \beta =(q+1) \alpha \beta$ and the injectivity of $\xi_{p,q}^{r}$ on $[p+1,q+1]$, we obtain $(q+1)\alpha \in \{p,q+1\}$. If $(q+1)\alpha = p$, then $\alpha_{\mid_{[p,q+1]}} = \xi_{p,q}^{r}$, which contradicts $\alpha \notin F^{r}_{p,q}$. If $(q+1) \alpha =q+1$, then $(q+1)\beta = p$, and so $\beta_{\mid_{[p,q+1]}} = \xi_{p,q}^{r}$, which contradicts $\beta \notin F^{r}_{p,q}$.

Let $\alpha \in F^{r}_{p,q}$. Since $\alpha_{\mid_{[p,q+1]}} = \xi_{p,q}^{r}$, we have $\xi_{p,q}^{r} = 1_{[p,q+1]}\alpha \in \langle (\mathcal{PORD}(n,r) \setminus F^{r}_{p,q}) \cup \{\alpha\}\rangle$, and so $E_{r}\cup F_{r}\cup G_{n} \subseteq \langle (\mathcal{PORD}(n,r) \setminus F^{r}_{p,q}) \cup \{\alpha\}\rangle$. Therefore, $\mathcal{PORD}(n,r) \setminus F^{r}_{p,q}$ is a maximal subsemigroup of $\mathcal{PORD}(n,r)$ since $\langle E_{r}\cup F_{r}\cup G_{n}\rangle = \mathcal{PORD}(n,r)$ by Theorem \ref{t8}. \qed

\begin{lemma}\label{l13}
	Let $n-\lfloor\frac{n}{3}\rfloor \leq r \leq n-1$ and $1\leq p\leq q-2 \leq n-2$ with $(p,q) \neq (1,n)$. Then $\mathcal{PORD}(n,r)\setminus G_{p,q}$ is a maximal subsemigroup of $\mathcal{PORD}(n,r)$.
\end{lemma}

\proof Assume that there are $\alpha, \beta\in \mathcal{PORD}(n,r) \setminus G_{p,q}$ such that $\alpha \beta\in G_{p,q}$. Then $\alpha \beta_{\mid_{\dom(\gamma_{p,q})}} = \gamma_{p,q}$, i.e. $1_{\dom(\gamma_{p,q})} \alpha \beta= \gamma_{p,q}$, and so by Lemma \ref{l2}, we have $x\alpha=x$ and $x\beta=x$ for all $x\in \{p,q\}$ since $\fix(\gamma_{p,q}) = \{p,q\}$. Then, since $\dom(\gamma_{p,q}) \subseteq \dom(\alpha)$, we have
$$\alpha_{\mid_{[p,p+k]}} = \left(
\begin{array}{cccc}
	p & p+1 & \cdots & p+k \\
	p & p-1 & \cdots & p-k \\
\end{array}
\right) \quad \mbox{ or } \quad \alpha_{\mid_{[p,p+k]}} = 1_{[p,p+k]}.$$
On the other hand
$$\alpha_{\mid_{[q,q+l]}} = \left(
\begin{array}{cccc}
	q & q+1 & \cdots & q+l \\
	q & q-1 & \cdots & q-l \\
\end{array}
\right) \quad \mbox{ or } \quad \alpha_{\mid_{[q,q+l]}} = 1_{[q,q+l]}.$$
Since $\alpha \in \mathcal{PORD}(n,r)$, we can conclude that $\alpha_{\mid_{\dom(\gamma_{p,q})}}$ is the partial identity on $\dom(\gamma_{p,q})$ or $\alpha_{\mid_{\dom(\gamma_{p,q})}} = \gamma_{p,q}$. Since $\alpha \notin G_{p,q}$, we obtain $\alpha_{\mid_{\dom(\gamma_{p,q})}} = 1_{\dom(\gamma_{p,q})}$. Hence, $\gamma_{p,q} = \alpha\beta_{\mid_{\dom(\gamma_{p,q})}} = \beta_{\mid_{\dom(\gamma_{p,q})}}$, which contradicts $\beta \notin G_{p,q}$. Consequently, we obtain that $\mathcal{PORD}(n,r)\setminus G_{p,q}$ is a subsemigroup of $\mathcal{PORD}(n,r)$.

Now, let $\alpha \in G_{p,q}$. Since $\alpha_{\mid_{\dom(\gamma_{p,q})}} =\gamma_{p,q}$, we have $\gamma_{p,q} = 1_{\dom(\gamma_{p,q})}\alpha \in \langle (\mathcal{PORD}(n,r) \setminus G_{p,q}) \cup\{ \alpha\} \rangle$, and so $E_{r}\cup F_{r}\cup G_{n} \subseteq \langle (\mathcal{PORD}(n,r) \setminus G_{p,q}) \cup \{\alpha\}\rangle$. Therefore, $\mathcal{PORD}(n,r) \setminus G_{p,q}$ is a maximal subsemigroup of $\mathcal{PORD}(n,r)$ since $\langle E_{r}\cup F_{r}\cup G_{n}\rangle = \mathcal{PORD}(n,r)$ by Theorem \ref{t8}.
\qed

Similarly to Lemma \ref{l13}, one can prove the following.

\begin{lemma}\label{l14}
	Let $3\leq r<n- \lfloor \frac{n}{3} \rfloor$,\, $1\leq p\leq q-2 \leq n-2$, and $1\leq s\leq \min\{p, q-p, r-1\}$ with $(q,s) \neq (n,1)$. Then $\mathcal{PORD}(n,r)\setminus H_{p,q}^{r,s}$ is a maximal subsemigroup of $\mathcal{PORD}(n,r)$. \hfill $\square$
\end{lemma}

Now, we can characterize the maximal subsemigroups of $\mathcal{PORD}(n,r)$ for $3 \leq r \leq n$.

\begin{theorem}\label{t15}
	Let $n-\lfloor\frac{n}{3}\rfloor \leq r \leq n-1$ and let $S$ be a subsemigroup of $\mathcal{PORD}(n,r)$. Then $S$ is a maximal subsemigroup of $\mathcal{PORD}(n,r)$ if and only if $S$ belongs to one of the following types:
	\begin{enumerate}
		\item[$(1)$] $S= \mathcal{PORD}(n,r) \setminus \{\varepsilon\}$, for each $\varepsilon \in E_{r}$.
		\item[$(2)$] $S= \mathcal{PORD}(n,r) \setminus F^{r}_{p,q}$, for all $1\leq p\leq n-2$ and $p+1\leq q\leq p+r-2\leq n-1$.
		\item[$(3)$] $S= \mathcal{PORD}(n,r) \setminus G_{p,q}$, for all $1\leq p\leq q-2 \leq n-2$ with $(p,q) \neq (1,n)$.
	\end{enumerate}
\end{theorem}

\proof Let $\varepsilon \in E_{r}$. By Proposition \ref{p4}, since $\varepsilon$ is undecomposable in $\mathcal{PORD}(n,r)$, we can conclude that $\mathcal{PORD}(n,r) \setminus\{\varepsilon\}$ is a maximal subsemigroup of $\mathcal{PORD}(n,r)$. By Lemma \ref{l12}, we know that $\mathcal{PORD}(n,r)\setminus F^{r}_{p,q}$ is a maximal subsemigroup of $\mathcal{PORD}(n,r)$ for all $1\leq p\leq n-2$ and $p+1\leq q\leq p+r-2\leq n-1$. By Lemma \ref{l13}, we know that $\mathcal{PORD}(n,r)\setminus G_{p,q}$ is a maximal subsemigroup of $\mathcal{PORD}(n,r)$ for all $1\leq p\leq q-2 \leq n-2$.

Conversely, let $S$ be a maximal subsemigroup of $\mathcal{PORD}(n,r)$. Suppose that $E_{r}\nsubseteq S$, i.e. there is $\varepsilon \in E_{r} \setminus S$. Since $\mathcal{PORD}(n,r) \setminus \{\varepsilon \}$ is a maximal subsemigroup of $\mathcal{PORD}(n,r)$, we obtain that $S=\mathcal{PORD}(n,r) \setminus \{\varepsilon \}$.

Suppose that $E_{r}\subseteq S$ and $F^{r}_{p,q}\cap S = \emptyset$ for some $1\leq p\leq n-2$ and $p+1\leq q\leq p+r-2\leq n-1$. Since $\mathcal{PORD}(n,r) \setminus F^{r}_{p,q}$ is a maximal subsemigroup of $\mathcal{PORD}(n,r)$, we obtain that $S = \mathcal{PORD}(n,r)\setminus F^{r}_{p,q}$.

Suppose that $E_{r}\subseteq S$ and $F^{r}_{p,q} \cap S \neq \emptyset$, i.e. $S$ contains at least one $\alpha \in \mathcal{PORD}(n,r)$ such that $\alpha_{\mid_{[p,q+1]}} = \xi_{p,q}^{r}$, for all $1\leq p\leq n-2$ and $p+1\leq q\leq p+r-2\leq n-1$. Thus, since $1_{[p,q+1]} \in \mathcal{PC}(n,r)=\langle E_{r}\rangle \subseteq S$, we have $1_{[p,q+1]}\alpha = \xi_{p,q}^{r} \in S$ for all $1\leq p\leq n-2$ and $p+1\leq q\leq p+r-2\leq n-1$. Therefore, we obtain $F_{r} \subseteq S$. Since $\mathcal{PORD}(n,r) = \langle E_{r}\cup F_{r}\cup G_{n}\rangle$ and $E_{r}\cup F_{r} \subseteq S$, we have that $G_{n} \nsubseteq S$.

Assume that $G_{p,q} \cap S \neq \emptyset$, i.e. $S$ contains at least one $\alpha \in \mathcal{PORD}^*(n,r)$ such that $\alpha_{\mid_{\dom(\gamma_{p,q})}} =\gamma_{p,q}$, for all $1\leq p\leq q-2 \leq n-2$ with $(p,q) \neq (1,n)$. Thus, since $1_{\dom(\gamma_{p,q})} \in S$, we have $\gamma_{p,q} = 1_{\dom(\gamma_{p,q})}\alpha \in S$ for all $1\leq p\leq q-2 \leq n-2$ with $(p,q) \neq (1,n)$. Therefore, we obtain $G_{n} \subseteq S$, which is a contradiction. Hence, there are $1\leq p\leq q-2 \leq n-2$ with $(p,q) \neq (1,n)$ such that $G_{p,q}\cap S= \emptyset$. Finally, we obtain $S = \mathcal{PORD}(n,r)\setminus G_{p,q}$, since $\mathcal{PORD}(n,r) \setminus G_{p,q}$ is a maximal subsemigroup of $\mathcal{PORD}(n,r)$.
\qed

Clearly, $\mathcal{PORD}_n = \mathcal{PORD}(n,n-1) \cup \{1_n\}$. Hence, we obtain the following result:

\begin{theorem}\label{t16}
	Let $S$ be a subsemigroup of $\mathcal{PORD}_n$. Then $S$ is a maximal subsemigroup of $\mathcal{PORD}_n$ if and only if $S = \mathcal{PORD}(n,n-1)$ or there exists a maximal subsemigroup $T$ of $\mathcal{PORD}(n,n-1)$ such that $S = T \cup \{1_n\}$. \hfill $\square$
\end{theorem}

Similarly to Theorem \ref{t15}, one can prove the following theorem.

\begin{theorem}\label{t17}
	Let $3 \leq r < n-\lfloor\frac{n}{3}\rfloor$ and let $S$ be a subsemigroup of $\mathcal{PORD}(n,r)$. Then $S$ is a maximal subsemigroup of $\mathcal{PORD}(n,r)$ if and only if $S$ belongs to one of the following types:
	\begin{enumerate}
		\item[$(1)$] $S= \mathcal{PORD}(n,r) \setminus \{\varepsilon\}$, for each $\varepsilon \in E_{r}$.
		\item[$(2)$] $S= \mathcal{PORD}(n,r) \setminus F^{r}_{p,q}$, for all $1\leq p\leq n-2$ and $p+1\leq q\leq p+r-2\leq n-1$.
		\item[$(3)$] $S= \mathcal{PORD}(n,r) \setminus H_{p,q}^{r,s}$, for all $1\leq p\leq q-2 \leq n-2$ and $1\leq s\leq \min\{p, q-p, r-1\}$ with $(q,s) \neq (n,1)$. \hfill $\square$
	\end{enumerate}
\end{theorem}

Now, we will characterize the maximal subsemigroups of $\mathcal{IORD}(n,r)$.

First, let $n-\lfloor\frac{n}{3}\rfloor \leq r \leq n-1$.
Since $EI_{r}\cup FI_{r}\cup GI_{r}\cup \left( \bigcup\limits_{k=2}^{r-1}GI^{c}_{k}\right) \cup G_{n}$ is a generating set of $\mathcal{IORD}(n,r)$ consisting entirely of undecomposable elements, by Lemmas \ref{il2} and \ref{il4}, we get immediately the following:
\begin{theorem}\label{it8}
Let $n-\lfloor\frac{n}{3}\rfloor \leq r \leq n-1$ and let $S$ be a subsemigroup of $\mathcal{IORD}(n,r)$. Then $S$ is a maximal subsemigroup of $\mathcal{IORD}(n,r)$ if and only if $S = \mathcal{IORD}(n,r)\setminus \{\alpha\}$, for each $\alpha \in EI_{r}\cup FI_{r}\cup GI_{r}\cup \left( \bigcup\limits_{k=2}^{r-1}GI^{c}_{k}\right) \cup G_{n}$. \hfill $\square$
\end{theorem}

Clearly, $\mathcal{IORD}_n = \mathcal{IORD}(n,n-1) \cup \{1_n\}$. Hence, we obtain the following result:
\begin{theorem}\label{it9}
Let $S$ be a subsemigroup of $\mathcal{IORD}_n$. Then $S$ is a maximal subsemigroup of $\mathcal{IORD}_n$ if and only if $S = \mathcal{IORD}_n \setminus \{\alpha\}$ for each $\alpha \in EI_{n-1}\cup FI_{n-1}\cup GI_{n-1}\cup \left( \bigcup\limits_{k=2}^{n-2}GI^{c}_{k}\right) \cup G_{n} \cup \{1_n\}$. \hfill $\square$
\end{theorem}

Now, let $3 \leq r < n-\lfloor\frac{n}{3}\rfloor$. For $1 \leq p \leq q-2 \leq n-2$ and $1 \leq s \leq \min\{p, q-p, r-1\}$ with $(q,s) \neq (n,1)$, let
$$HI_{p,q}^{r,s} = H_{p,q}^{r,s} \cap \mathcal{IORD}(n,r).$$

\begin{lemma}\label{il10}
Let $3 \leq r < n-\lfloor\frac{n}{3}\rfloor$, $1 \leq p \leq q-2 \leq n-2$, and $1 \leq s \leq \min\{p, q-p, r-1\}$ with $(q,s) \neq (n,1)$. Then $\mathcal{IORD}(n,r)\setminus HI_{p,q}^{r,s}$ is a maximal subsemigroup of $\mathcal{IORD}(n,r)$.
\end{lemma}

\proof
Let $Y = \dom(\gamma_{p,q}^{r,s}) = [p,p+s-1]\cup[q,q+u-1]$, where $u=\min\{r-s,q-p,n-q+1\}$. Obviously, we have $HI_{p,q}^{r,s} = \{\gamma_{p,q}^{r,s}\}$ if $\min\{r-s, q-p, n-q+1\} = r-s$, i.e. $|Y|=r$. So in this case, $\mathcal{IORD}(n,r)\setminus HI_{p,q}^{r,s} = \mathcal{IORD}(n,r)\setminus \{\gamma_{p,q}^{r,s}\}$ is a maximal subsemigroup of $\mathcal{IORD}(n,r)$, since $\gamma_{p,q}^{r,s}$ is undecomposable in $\mathcal{IORD}(n,r)$ by Lemma \ref{il7}.

Let now $\lvert Y\rvert<r$. Assume that there are $\alpha, \beta\in \mathcal{IORD}(n,r) \setminus HI_{p,q}^{r,s}$ such that $\alpha \beta\in HI_{p,q}^{r,s}$. Then $\alpha\beta_{\mid_{Y}} = \gamma_{p,q}^{r,s}$, i.e. $1_{Y} \alpha \beta = \gamma_{p,q}^{r,s}$, and so by Lemma \ref{l2}, we have $x\alpha=x$ and $x\beta=x$ for all $x\in \{p,q\}$ since $\fix(\gamma_{p,q}^{r,s}) = \{p,q\}$. Then, since $\dom(\gamma_{p,q}^{r,s}) = [p,p+s-1]\cup[q,q+u-1] \subseteq \dom(\alpha)$, we have
$$\alpha_{\mid_{[p,p+s-1]}} = \left(
\begin{array}{cccc}
	p & p+1 & \cdots & p+s-1 \\
	p & p-1 & \cdots & p-s+1 \\
\end{array}
\right) \quad \mbox{ or } \quad \alpha_{\mid_{[p,p+s-1]}} = 1_{[p,p+s-1]}.$$
On the other hand
$$\alpha_{\mid_{[q,q+u-1]}} = \left(
\begin{array}{cccc}
	q & q+1 & \cdots & q+u-1 \\
	q & q-1 & \cdots & q-u+1 \\
\end{array}
\right) \quad \mbox{ or } \quad \alpha_{\mid_{[q,q+u-1]}} = 1_{[q,q+u-1]}.$$
Since $\alpha \in \mathcal{IORD}(n,r)$, we can conclude that $\alpha_{\mid_{Y}}$ is the partial identity on $\dom(\gamma_{p,q}^{r,s})$ or $\alpha_{\mid_{Y}} = \gamma_{p,q}^{r,s}$. Since $\alpha \notin HI_{p,q}^{r,s}$, we obtain $\alpha_{\mid_{Y}} = 1_{\dom(\gamma_{p,q}^{r,s})}$. Hence, $\gamma_{p,q}^{r,s} = \alpha\beta_{\mid_{Y}} = \beta_{\mid_{\dom(\gamma_{p,q}^{r,s})}} = \gamma_{p,q}^{r,s}$,  which contradicts $\beta \notin HI_{p,q}^{r,s}$.
Consequently, we obtain that $\mathcal{IORD}(n,r)\setminus HI_{p,q}^{r,s}$ is a subsemigroup of $\mathcal{IORD}(n,r)$.

Let $\alpha \in HI_{p,q}^{r,s}$. Since $\alpha_{\mid_{Y}} = \gamma_{p,q}^{r,s}$, we have $\gamma_{p,q}^{r,s} = 1_{Y}\alpha \in \langle (\mathcal{IORD}(n,r) \setminus HI_{p,q}^{r,s}) \cup\{ \alpha\} \rangle$, and so $EI_{r}\cup FI_{r}\cup GI_{r}\cup \left( \bigcup\limits_{k=2}^{r-1}GI^{c}_{k}\right)\cup H_{n}^{r} \subseteq \langle (\mathcal{IORD}(n,r) \setminus HI_{p,q}^{r,s}) \cup \{\alpha\}\rangle$. Therefore, $\mathcal{IORD}(n,r) \setminus HI_{p,q}^{r,s}$ is a maximal subsemigroup of $\mathcal{IORD}(n,r)$ since $\langle EI_{r}\cup FI_{r}\cup GI_{r}\cup \left( \bigcup\limits_{k=2}^{r-1}GI^{c}_{k}\right)\cup H_{n}^{r}\rangle = \mathcal{IORD}(n,r)$ by Theorem \ref{it6}.
\qed

\begin{theorem}\label{it11}
Let $3 \leq r < n-\lfloor\frac{n}{3}\rfloor$ and let $S$ be a subsemigroup of $\mathcal{IORD}(n,r)$. Then $S$ is a maximal subsemigroup of $\mathcal{IORD}(n,r)$ if and only if $S$ belongs to one of the following types:
\begin{enumerate}
  \item[(1)] $S = \mathcal{IORD}(n,r)\setminus \{\alpha\}$, for each $\alpha \in EI_{r}\cup FI_{r}\cup GI_{r}\cup \left( \bigcup\limits_{k=2}^{r-1}GI^{c}_{k}\right)$.\\
  \item[(2)] $S = \mathcal{IORD}(n,r)\setminus HI_{p,q}^{r,s}$, for all $1 \leq p \leq q-2 \leq n-2$ and $1 \leq s \leq \min\{p, q-p, r-1\}$ with $(q,s) \neq (n,1)$.\\
\end{enumerate}
\end{theorem}

\proof
Let $\alpha \in EI_{r}\cup FI_{r}\cup GI_{r}\cup \left( \bigcup\limits_{k=2}^{r-1}GI^{c}_{k}\right)$. Since $\alpha$ is undecomposable in $\mathcal{IORD}(n,r)$ by Lemma \ref{il2}, we can conclude that $\mathcal{IORD}(n,r) \setminus\{\alpha\}$ is a maximal subsemigroup of $\mathcal{IORD}(n,r)$. By Lemma \ref{il10}, we know that $\mathcal{IORD}(n,r)\setminus HI_{p,q}^{r,s}$ is a maximal subsemigroup of $\mathcal{IORD}(n,r)$, for all $1 \leq p \leq q-2 \leq n-2$ and $1 \leq s \leq \min\{p, q-p, r-1\}$ with $(q,s) \neq (n,1)$.

Conversely, let $S$ be a maximal subsemigroup of $\mathcal{IORD}(n,r)$. Suppose that $(EI_{r}\cup FI_{r}\cup GI_{r}\cup \left( \bigcup\limits_{k=2}^{r-1}GI^{c}_{k}\right)) \nsubseteq S$, i.e. there is $\alpha \in (EI_{r}\cup FI_{r}\cup GI_{r}\cup \left( \bigcup\limits_{k=2}^{r-1}GI^{c}_{k}\right)) \setminus S$. Since $\mathcal{IORD}(n,r) \setminus \{\alpha \}$ is a maximal subsemigroup of $\mathcal{IORD}(n,r)$, we obtain that $S=\mathcal{IORD}(n,r) \setminus \{\alpha \}$.

Suppose that $(EI_{r}\cup FI_{r}\cup GI_{r}\cup \left( \bigcup\limits_{k=2}^{r-1}GI^{c}_{k}\right))\subseteq S$ and $HI_{p,q}^{r,s} \cap S \neq \emptyset$, i.e. $S$ contains at least one $\alpha \in \mathcal{IORD}^*(n,r)$ such that $\alpha_{\mid_{\dom(\gamma_{p,q}^{r,s})}} = \gamma_{p,q}^{r,s}$, for all $1 \leq p \leq q-2 \leq n-2$ and $1 \leq s \leq \min\{p, q-p, r-1\}$ with $(q,s) \neq (n,1)$. Thus, since $1_{\dom(\gamma_{p,q}^{r,s})} \in IC(n,r) = \langle EI_{r} \cup FI_{r}\rangle \subseteq S$, we have $\gamma_{p,q}^{r,s} = 1_{\dom(\gamma_{p,q}^{r,s})}\alpha \in S$, for all $1 \leq p \leq q-2 \leq n-2$ and $1 \leq s \leq \min\{p, q-p, r-1\}$ with $(q,s) \neq (n,1)$. Therefore, we obtain $H_{n}^{r} \subseteq S$. Since $\mathcal{IORD}(n,r) = \langle EI_{r}\cup FI_{r}\cup GI_{r}\cup \left( \bigcup\limits_{k=2}^{r-1}GI^{c}_{k}\right)\cup H_{n}^{r}\rangle$ and $(EI_{r}\cup FI_{r}\cup GI_{r}\cup \left( \bigcup\limits_{k=2}^{r-1}GI^{c}_{k}\right)) \subseteq S$, we obtain a contradiction with maximality of $S$.
Hence, there are $1 \leq p \leq q-2 \leq n-2$ and $1 \leq s \leq \min\{p, q-p, r-1\}$  with $(q,s) \neq (n,1)$ such that $HI_{p,q}^{r,s} \cap S= \emptyset$. Finally, we obtain $S = \mathcal{IORD}(n,r)\setminus HI_{p,q}^{r,s}$, since $\mathcal{IORD}(n,r) \setminus HI_{p,q}^{r,s}$ is a maximal subsemigroup of $\mathcal{IORD}(n,r)$.
\qed

Since we have not found any explicit formula for the cardinality of $\lvert H_{n}^{r}\rvert$, we have not written the ranks of $\mathcal{PORD}(n,r)$ and $\mathcal{IORD}(n,r)$ for the case $3\leq r< n-\lfloor\frac{n}{3}\rfloor$.

\vspace{3mm}
\noindent
\textbf{Open Problem.} Does there exist an explicit formula for the cardinality of $\lvert H_{n}^{r}\rvert$ for each $3\leq r< n-\lfloor\frac{n}{3}\rfloor$?

\section*{Declaration of Interest}

No potential conflict of interest was reported by the authors.

\end{document}